\pdfoutput=1
\documentclass{article}
\usepackage[left=1in,right=1in,top=1in,bottom=1.3in]{geometry}
\usepackage{amsmath, amssymb, amsthm, latexsym, amsfonts, indentfirst, xcolor, mathtools, manfnt, microtype, enumitem, graphicx}
\usepackage[utf8]{inputenc}
\usepackage[english]{babel}
\usepackage{csquotes}
\graphicspath{ {./} }

\usepackage[breaklinks]{hyperref}
\usepackage{cleveref}
\hypersetup{
	colorlinks = true, % Colours links instead of ugly boxes
	urlcolor = cyan, % Colour for external hyperlinks
	linkcolor = teal, % Colour of internal links
	citecolor = cyan % Colour of citations
}

\usepackage[backend=biber, sorting=nty, giveninits=true, maxbibnames=99]{biblatex}
\renewbibmacro{in:}{}
\DeclareFieldFormat{pages}{#1}

\usepackage{subfiles}

\usepackage{enumitem}
\newlist{property}{enumerate}{1}
\setlist[property,1]{label=(P\arabic*), ref=P\arabic*}
\newlist{observation}{enumerate}{1}
\setlist[observation,1]{label=(O\arabic*), ref=O\arabic*}

\newtheorem{Conjecture}{Conjecture} 
\newtheorem{Corollary}{Corollary} 
 
\newtheorem{Lemma}{Lemma} 
\newtheorem{Proposition}{Proposition} 
\newtheorem{Theorem}{Theorem}

\newtheorem{Definition}{Definition}

\newcommand{\Z}{{\mathbb Z}}
\newcommand{\N}{{\mathbb N}}

\newcommand{\sm}{\!\setminus\!}

\newcommand{\upperAlpha}{{\alpha}}
\newcommand{\upperBeta}{{\beta}}

\newcommand{\upperXi}{{i}}
\newcommand{\upperEta}{{j}}
\newcommand{\upperTau}{{\upperEta}}

\newcommand{\upperAA}{{A}}
\newcommand{\upperA}{{b}}
\newcommand{\upperB}{{a}}
\newcommand{\upperCC}{{C}}

\newcommand{\upperD}{{d}}

\newcommand{\upperI}{{n}}
\newcommand{\upperJ}{{t}}
\newcommand{\upperII}{{I}}
\newcommand{\upperK}{{m}}
\newcommand{\upperMM}{{M}}
\newcommand{\upperM}{{k}}
\newcommand{\upperNN}{{N}}
\newcommand{\upperN}{{n}}
\newcommand{\upperSS}{{S}}

\newcommand{\upperR}{{r}}
\newcommand{\upperTT}{{T}}

\newcommand{\upperUU}{{U}}

\newcommand{\upperVV}{{V}}
\newcommand{\upperV}{{v}}
\newcommand{\upperWW}{{W}}
\newcommand{\upperW}{{w}}

\newcommand{\lowera}{a}
\newcommand{\lowerb}{b}
\newcommand{\lowerk}{k}
\newcommand{\lowerm}{m}
\newcommand{\loweri}{i}
\newcommand{\lowerj}{j}
\newcommand{\lowerd}{d}
\newcommand{\lowerr}{r}
\newcommand{\lowern}{n}
\newcommand{\lowerz}{z}
\newcommand{\lowert}{t}
\newcommand{\lowerp}{p}

\newcommand{\theora}{a}
\newcommand{\theorb}{b}
\newcommand{\theork}{k}
\newcommand{\theorm}{m}

\title{Packing Density of Sets With Only Two Nonmixed Gaps}
\author{
    Alexander Natalchenko\thanks{MIPT, Moscow, Russia.  Email:~\href{mailto:natalchenko.ae@gmail.com}{\tt natalchenko.ae@gmail.com}}
    \and
    Arsenii Sagdeev\thanks{Alfréd Rényi Institute of Mathematics, Budapest, Hungary.  Email:~\href{mailto:sagdeevarsenii@gmail.com}{\tt sagdeevarsenii@gmail.com}.}
}

\date{}

\begin{document}

\maketitle

\begin{abstract}
 For a finite set of integers such that the first few gaps between its consecutive elements equal $a$, while the remaining gaps equal $b$, we study dense packings of its translates on the line. We obtain an explicit lower bound on the corresponding optimal density, conjecture its tightness, and prove it in case one of the gap lengths, $a$ or $b$, appears only once. This is equivalent to a Motzkin problem on the independence ratio of certain integer distance graphs.
\end{abstract}

\section{Introduction}

Introduced by Arthur Cayley in 1878, \textit{Cayley graphs} play an important role in modern mathematics as a link between its different branches such as Algebra, Combinatorics, and Number Theory~\cite{alon2013chromatic, babai1979spectra,cervantes2023chromatic,konstantinova2008some, li2002isomorphisms}. In the present paper, we focus on the Cayley graphs corresponding to finite sets of positive integers. More specifically, for a finite $M \subset \N$, its \textit{Cayley graph} $G(\Z, M)$, which is also often called the \textit{integer distance graph}, is a graph with the vertex set $\Z$ where two vertices $v_1,v_2 \in \Z$ are adjacent if and only if $|v_1-v_2| \in M$. Plenty of research papers are dedicated to the study of various natural parameters of these graphs, foremost, of their chromatic numbers, see \cite{barajas2008distance,eggleton1985colouring,heuberger2003planarity,katznelson2001chromatic,oravcova2024problem,sokolov2023chromatic,zhu2002circular}.

Perhaps the second most studied parameter of these graphs is their \textit{independence ratio}, defined as the maximum upper density\footnote{Recall that the \textit{upper density} of $A \subset \Z$ is defined  as $\overline{d}(A)=\limsup_{n\to \infty}\frac{|\{a\in A: -n \le a \le n\}|}{2n+1}$.} of their independent subsets. Note that in the heart of the recent breakthroughs in combinatorial geometry due to Davies~\cite{davies2023chromatic, davies2024odd} and  Davies, McCarty, and Pilipczuk ~\cite{davies2023prime}, there is a resourceful \textit{estimation} of the independence ratios for some special families of Cayley graphs. This illustrates that the \textit{exact} determination of the independence ratio is very hard in general. A systematic study of this problem goes back to the paper~\cite{cantor1973sequences} of Cantor and Gordon from 1973, in which the authors attributed it to Motzkin and coined out the notation $\mu(M)$ for the independence ratio of $G(\Z,M)$. In other words, they defined $\mu(M)$ as the maximum upper density of an \textit{$M$-avoiding set} $A \subset \Z$, i.e., such that  $a_1 - a_2 \notin M$ for all $a_1,a_2 \in A$.

It is easy to see that $\mu(M) = \frac{1}{2}$ if $M$ is a singleton. For a two-point $M=\{a,b\}$ with coprimes $a$ and $b$, Cantor and Gordon~\cite{cantor1973sequences} showed %\footnote{We provide a short proof of this fact here for completeness. First, note that the set $C=\{0,a,2a,\dots,ba, b(a-1),\dots,b\}$ forms a cycle of length $a+b$ in $G(\Z, M)$, and thus $|A\cap C| \le \lfloor(a+b)/2\rfloor$ for every $M$-avoiding set $A$. By averaging over all translates of $C$, we conclude that $\mu(M) \le \frac{\lfloor(a+b)/2\rfloor}{a+b}$. For the matching lower bound, it is sufficient check that the set $A = \{2\lowert a: 0 \le \lowert < \lfloor(a+b)/2\rfloor\}+(a+b)\Z$ is $M$-avoiding and of density $\frac{\lfloor(a+b)/2\rfloor}{a+b}$.}
that $\mu(M) = \frac{\lfloor(a+b)/2\rfloor}{a+b}$. Despite the long history of research, already in case $|M|=3$, a closed-form expression for $\mu(M)$ was found only for some special families of triples, see~\cite{liu2020sequences} and the references therein. However, several such expressions were found for sets with additional structure, e.g., whose elements form an arithmetic~\cite{gupta2000sets} or a geometric~\cite{pandey2015note} progression or a union of two intervals~\cite{liu2008fractional}. Another notable example of this sort is $M=\{a,b,a+b\}$, which was initially studied by Rabinowitz and Proulx~\cite[Theorem~5.4]{rabinowitz1985asymptotic}, who found a lower bound on the corresponding $\mu(M)$ and conjectured its tightness. Their conjecture was later confirmed by Liu and Zhu in~\cite[Theorem~3.1]{liuzhu2004} and~\cite{liu2005d}.

\begin{Theorem} \label{th:liuzhu}
	Let $a,b$ be coprime positive integers. Then for $M=\{a,b,a+b\}$, we have
	\begin{equation*}
		\mu(M) = \max \left (\frac{\lfloor (a+2b)/3\rfloor}{a+2b}, \frac{\lfloor (2a+b)/3\rfloor}{2a+b}\right) =
		\begin{dcases}
			\frac{1}{3}, & a-b=3d, \\
			\frac{b+d}{a+2b}, & a-b=3d+1, \\
			\frac{a-d-1}{2a+b}, & a-b=3d+2.
		\end{dcases}
	\end{equation*}
\end{Theorem}

In the present paper, we generalize these results to $M=\{\loweri\lowera+\lowerj\lowerb: 0 \le \loweri \le \lowerk,\, 0 \le \lowerj \le \lowerm,\, \loweri+\lowerj>0\}$. Before diving into details, we present an alternative point of view on this problem, which was also a key motivation for our study.

For a finite $S \subset \Z$, a subset $A \subset \Z$ is called \textit{$S$-packing} if and only if two translates $a_1+S$ and $a_2+S$ of $S$ are disjoint for all distinct $a_1,a_2 \in A$. Naturally, the \textit{packing density} $d_p(S)$ of $S$ is defined as the maximum upper density of an $S$-packing set. The study of this parameter, along with its covering counterpart, goes back to Weinstein~\cite{weinstein1976some} and Newman~\cite{newman1967complements}, respectively. However, it is not hard to see\footnote{Indeed, two distinct translates $a_1+S$ and $a_2+S$ share a common point if and only if $a_1+s_1=a_2+s_2$ for some $s_1 \neq s_2 \in S$, which can be rewritten as $a_1-a_2 = s_2-s_1$.} that $A$ is $S$-packing if and only if $A$ is $M$-avoiding for $M=\{s_2-s_1: s_1,s_2 \in S, s_1<s_2\}$, and thus $d_p(S)=\mu(M)$. In particular, we have $d_p(S)=\mu(M)$ for $S=\{0,a,a+b\}$ and $M=\{a,b,a+b\}$. In this notation, the result of \Cref{th:liuzhu} was independently conjectured by Schmidt and Tuller~\cite{schmidt2008covering} and later reproved in~\cite{frankl2023solution}. 

In the present paper, we continue this line of research, and study the packing density of sets such that the first few gaps between their consecutive elements equal $a$, while the remaining gaps equal $b$. More specifically, we consider $S = \{0, \theora, \dots, \theork\theora, \theork\theora + \theorb, \dots, \theork\theora + \theorm\theorb\}$. As noted above, the equality $d_p(S) = \mu(M)$ holds for $M=\{\loweri\lowera+\lowerj\lowerb: 0 \le \loweri \le \lowerk,\, 0 \le \lowerj \le \lowerm,\, \loweri+\lowerj>0\}$. Our first result provides a lower bound on these densities. 

\begin{Theorem}\label{main-conjecture}
    Let $\theora, \theorb$ be coprime positive integers. For $\theork, \theorm \in \N$, let $d \in \Z$ and $0 \le r \le k+m$ be unique integers such that $\theora - \theorb = (\theork + \theorm + 1) d + r$. Then for $S = \{0, \theora, \dots, \theork\theora, \theork\theora + \theorb, \dots, \theork\theora + \theorm\theorb\}$ and $M=\{\loweri\lowera+\lowerj\lowerb: 0 \le \loweri \le \lowerk,\, 0 \le \lowerj \le \lowerm,\, \loweri+\lowerj>0\}$, we have
    \begin{equation*}
        d_p (S) = \mu(M) \ge
        \begin{dcases}
        	\frac{1}{\theork + \theorm + 1}, & r=0, \\
        	\frac{\theorb + \theork d}{\theork\theora+ (\theorm+1)\theorb}, & 1 \le r \le \theorm, \\
            \frac{\theora - \theorm (d+1)}{(\theork + 1)\theora+ \theorm\theorb}, & \theorm+1 \le r \le \theork + \theorm.
        \end{dcases}
    \end{equation*}
\end{Theorem}

Let us make a couple of clarifications about this statement. First, the condition that $a$ and $b$ are coprime does not detract from the generality because simultaneously dividing all elements of $S$ or $M$ by their greatest common divisor does not affect the corresponding densities. Second, in case $r=0$, the tightness of our lower bound follows from a trivial inequality $d_p(S) \le 1/|S|$. The other two cases are symmetric: simultaneously switching the roles of $a$ and $b$, as well as of $k$ and $m$, will also interchange the fractions because $b-a = -(k+m+1)(d+1)+(k+m+1-r)$. Such symmetry is not surprising since this switching does not change the set $M$ at all and replaces $S$ with its reflection, which clearly does not affects the packing density. Finally, note that, as earlier in \Cref{th:liuzhu}, the denominators of these fractions equal the maximum element of $M$, namely $ka+mb$, increased by $a$ or $b$ .

We conjecture that \Cref{main-conjecture} is always tight. The case $k=m=1$ is simply the statement of \Cref{th:liuzhu}. In the main result of the present paper, we confirm this conjecture if either $k$ or $m$ equals $1$. Due to a symmetry discussed above, we can assume without loss of generality that $m=1$.

\begin{Theorem}\label{main-theorem}
    Let $\theora, \theorb$ be coprime positive integers. For $\theork \in \N$, let $d \in \Z$ and $0 \le r \le k+1$ be unique integers such that $\theora - \theorb = (\theork + 2) d + r$. Then for $S = \{0, \theora, \dots, \theork\theora, \theork\theora + \theorb\}$ and $M=\{\theora, \dots, \theork\theora\}\cup \{\theorb, \theora+\theorb, \dots, \theork\theora + \theorb\}$, we have
    \begin{equation*}
    	d_p (S) = \mu(M) =
    	\begin{dcases}
    		\frac{1}{\theork +2}, & r=0, \\
    		\frac{\theorb + \theork d}{\theork\theora+ 2\theorb}, & r = 1, \\
    		\frac{\theora -d-1}{(\theork + 1)\theora+ \theorb}, & 2\le r \le \theork+1.
    	\end{dcases}
    \end{equation*}
\end{Theorem}

\noindent
\textbf{Paper outline.} In the next two sections, we prove Theorems~\ref{main-conjecture} and~\ref{main-theorem}, respectively. In \Cref{SecConc}, we discuss covering counterparts of the aforementioned packing densities in more detail and state several open problems.

\section{Lower bound --- Proof of Theorem~\ref{main-conjecture}}

Our argument in this section generalizes the ideas for the special case $k=m=1$ from \cite[Proposition~5.3]{rabinowitz1985asymptotic}.

First, we consider the case $r=0$, that is, $a$ and $b$ have the same residue modulo $k+m+1$. Observe that this residue is relatively prime to $k+m+1$, since otherwise $a$ and $b$ would have a common divisor. Therefore, $S$ intersects each residue class modulo $k+m+1$ exactly once. Now it is easy to see that the set $(k+m+1)\Z$ containing only multiples of $(k+m+1)$ is $S$-packing and of density $\frac{1}{k+m+1}$, as desired.

Recall from the introduction that the other two cases are symmetric, and thus we can assume without loss of generality that $1 \le \lowerr\le \lowerm$. Our goal is to find an $M$-avoiding set of density $\frac{\lowerb+\lowerk\lowerd}{\lowern}$, where $\lowern=\lowerk\lowera+(\lowerm+1)\lowerb$. We claim that the set%\footnote{It is easy to check that if $\lowerr=0$, then $\gcd(\lowerb-\lowera, \lowern) = \gcd\big((\lowerk+\lowerm+1)\lowerd, (\lowerk+\lowerm+1)(\lowera+\lowerm\lowerd)\big) = (\lowerk+\lowerm+1)\gcd(\lowerd,\lowera) = \lowerk+\lowerm+1$, where the last equality is due to the fact that $\lowera$ and $\lowerb$ are coprime. Therefore, $A$ consists of all the integers divisible by $\lowerk+\lowerm+1$, i.e., it coincides with the construction we have already discussed in the introduction.}
\begin{equation*}
	A = \{\lowert(\lowera-\lowerb): 0 \le \lowert < \lowerb+\lowerk\lowerd\}+\lowern\Z
\end{equation*}
satisfies the desired conditions.

Assume the contrary, namely that either the density of $A$ is smaller than $\frac{\lowerb+\lowerk\lowerd}{\lowern}$, or that $A$ is not $M$-avoiding. In the former case, some two of its `residues' are equal modulo $\lowern$, that is
\begin{equation*}
	\lowert_1(\lowera-\lowerb) - \lowert_2(\lowera-\lowerb)=\lowern\lowerz,
\end{equation*}
for some $\lowerz \in \Z$ and distinct $0 \le \lowert_1,\lowert_2 < \lowerb+\lowerk\lowerd$. Similarly, in the latter case, the difference between two of the `residues' belongs to $M$ modulo $\lowern$, that is
\begin{equation*}
	\lowert_1(\lowera-\lowerb) - \lowert_2(\lowera-\lowerb)=\loweri\lowera+\lowerj\lowerb+\lowern\lowerz,
\end{equation*}
for some $\lowerz \in \Z$, $0 \le \loweri \le \lowerk,\, 0 \le \lowerj \le \lowerm$ such that $\loweri+\lowerj>0$, and not necessarily distinct $0 \le \lowert_1,\lowert_2 < \lowerb+\lowerk\lowerd$.
For both the former and the latter cases, we denote the difference $\lowert_1-\lowert_2$ by $\lowert$ and together we obtain
\begin{equation} \label{eq:low0}
	\lowert(\lowera-\lowerb)=\loweri\lowera+\lowerj\lowerb+\lowern\lowerz,
\end{equation}
where $-(\lowerb+\lowerk\lowerd) < \lowert < \lowerb+\lowerk\lowerd$, and at least one of the variables $\loweri,\lowerj,\lowerz$ does not equal $0$.

To get a desired contradiction, we first substitute the value of $\lowern=\lowerk\lowera+(\lowerm+1)\lowerb$ into \eqref{eq:low0}, and then collect the multiples of $\lowera$ and $\lowerb$ in different sides of the equation:
\begin{equation*}
	\lowera(\lowert-\lowerk\lowerz-\loweri)= \lowerb(\lowert+(\lowerm+1)\lowerz+\lowerj).
\end{equation*}
Since $\lowera$ and $\lowerb$ are coprime, there exists some $\lowerp \in \Z$ such that
\begin{equation} \label{eq:low1}
	\begin{cases}
		\lowera\lowerp = \lowert+(\lowerm+1)\lowerz+\lowerj, \\
		\lowerb\lowerp = \lowert-\lowerk\lowerz-\loweri.
	\end{cases}
\end{equation}
Subtracting the second equation from the first one yields that
\begin{equation} \label{eq:low2}
	\lowerp(\lowera-\lowerb)=(\lowerk+\lowerm+1)\lowerz+\loweri+\lowerj.
\end{equation} 
Recall that $\lowera-\lowerb=(\lowerk+\lowerm+1)\lowerd+\lowerr$ and put $\lowerr' = \loweri+\lowerj$ and $\lowerd' = \lowerz-\lowerd\lowerp$. In this notation, \eqref{eq:low2} is equivalent to
\begin{equation} \label{eq:low3}
	\lowerp\lowerr = (\lowerk+\lowerm+1)\lowerd'+\lowerr'.
\end{equation}

Next, we multiply the first equation of \eqref{eq:low1} by $\lowerk$, the second one by $\lowerm+1$, and sum these products:
\begin{equation} \label{eq:low4}
	\lowerp(\lowerk\lowera+(\lowerm+1)\lowerb)=(\lowerk+\lowerm+1)\lowert+\lowerk\lowerj-(\lowerm+1)\loweri.
\end{equation}
Then we rewrite the right-hand side of the last equation as
\begin{equation*}
	(\lowerk+\lowerm+1)\lowert+\lowerk\lowerj-(\lowerm+1)\loweri = (\lowerk+\lowerm+1)(\lowert-\loweri)+\lowerk\lowerr'
\end{equation*}
and use \eqref{eq:low3} to rewrite the left-hand side of \eqref{eq:low4} as
\begin{align*}
	\lowerp(\lowerk\lowera+(\lowerm+1)\lowerb) =&\ (\lowerk+\lowerm+1)\lowerb\lowerp+\lowerk\lowerp(\lowera-\lowerb) \\ 
	=&\ (\lowerk+\lowerm+1)\lowerb\lowerp+(\lowerk+\lowerm+1)\lowerk\lowerd\lowerp+\lowerk\lowerp\lowerr \\
	=&\ (\lowerk+\lowerm+1)(\lowerb\lowerp+\lowerk\lowerd\lowerp)+\lowerk((\lowerk+\lowerm+1)\lowerd'+\lowerr') \\
	=&\ (\lowerk+\lowerm+1)(\lowerb\lowerp+\lowerk\lowerd\lowerp+\lowerd'\lowerk)+\lowerk\lowerr'.
\end{align*}
After these straightforward rearrangements, we see that \eqref{eq:low4} is equivalent to
\begin{equation*}
	(\lowerk+\lowerm+1)(\lowerb\lowerp+\lowerk\lowerd\lowerp+\lowerd'\lowerk)+\lowerk\lowerr' = (\lowerk+\lowerm+1)(\lowert-\loweri)+\lowerk\lowerr',
\end{equation*}
and thus
\begin{equation} \label{eq:low5}
	\lowert=\lowerp(\lowerb+\lowerk\lowerd)+\lowerd'\lowerk+\loweri.
\end{equation}

Let us also analyze how the sign of $\lowerd'\lowerk+\loweri$ depends on the sign of $\lowerp$. Note that $0 \le \lowerr' = \loweri+\lowerj \le \lowerk+\lowerm$ by definition, and so \eqref{eq:low3} represents a standard division of $\lowerp\lowerr$ by $\lowerk+\lowerm+1$ with $\lowerd'$ being the quotient and $\lowerr'$ being the residue. In particular, if $\lowerp>0$, then $\lowerd' \ge 0$, and thus $\lowerd'\lowerk+\loweri \ge 0$. Similarly, if $\lowerp<0$, then $\lowerd'<0$, and thus $\lowerd'\lowerk+\loweri\le \loweri-\lowerk \le 0$.  %In the last remaining case, if $p = 0$, we obtain $d' = r' = 0.$

Finally, we consider three separate cases depending on the sign of $\lowerp$, and quickly get a desired contradiction in each of them. First, if $\lowerp>0$, then \eqref{eq:low5} implies that $\lowert \ge \lowerb+\lowerk\lowerd$, a contradiction with the inequalities on $\lowert$ posed in \eqref{eq:low0}. Similarly, if $\lowerp<0$, then \eqref{eq:low5} implies that $\lowert \le -(\lowerb+\lowerk\lowerd)$, a contradiction again. In the last remaining case, if $\lowerp=0$, we obtain $z = 0$ and $i + j = 0$ by~\eqref{eq:low2}. Therefore, $\loweri=\lowerj=\lowerz=0$, which contradicts the assumption that at least one of these variables does not equal $0$.

\section{Upper bound --- Proof of Theorem~\ref{main-theorem}}

Our argument in this section is built on the ideas of Liu and Zhu from~\cite[Section~4]{liuzhu2004} and~\cite{liu2005d} for the special case $k=m=1$, so we adhere to their notation whenever possible.

\subsection[The main inequality]{Reduction to an auxiliary inequality} \label{Sec3.1}

In this subsection, we work in the general setting of \Cref{main-conjecture}, where
\begin{equation*}
	S = \{0, \theora, \dots, \theork\theora, \theork\theora + \theorb, \dots, \theork\theora + \theorm\theorb\}, \ \ \ M=\{\loweri\lowera+\lowerj\lowerb: 0 \le \loweri \le \lowerk,\, 0 \le \lowerj \le \lowerm,\, \loweri+\lowerj>0\},
\end{equation*}
and show that the matching upper bound on $d_p(S)=\mu(M)$ follows from some auxiliary inequality. In the next two subsections, we prove this inequality in a special case either $k$ or $m$ equals $1$, thereby completing the proof of \Cref{main-theorem}.

As we discussed in the introduction, \Cref{main-conjecture} is symmetric, so in what follows, we can assume without loss of generality %\footnote{Note that \Cref{main-theorem} is \textit{not} symmetric in this sense, and we consider the other option $\upperB<\upperA$ in \Cref{Sec3.3}.}
that $\upperB>\upperA$. For $x,y \in \Z$, we denote an `integral segment' $\{z \in \Z: x \le z \le y\}$ by $[x,y]$, or simply by $[y]$ if $x=1$. For $\upperAlpha \in \Z$, let $\upperSS_\upperAlpha = \upperAlpha + \upperSS$ be a translate of $\upperSS$. 

\begin{Definition}\label{def_sets}
    Given an $\upperMM$-avoiding set $\upperAA \subset \mathbb{Z}$ such that $0\in \upperAA$, we also consider the following auxiliary sets:
    \begin{align*}
	& \upperII = \{ \upperAlpha \in [\upperA - 1] : \upperAA \cap \upperSS_\upperAlpha = \varnothing \}, \\
	& \upperTT_i = \upperAA \cap [i\upperB + \upperA + 1, (i + 1)\upperB - 1], \;\; 0 \le i \le \upperM - 1, \;\; T = \bigcup_{i = 0}^{\upperM - 1} \upperTT_i, \\
	& \upperUU = [\upperM \upperB + (\upperK + 1)\upperA, (\upperM + 1)\upperB + \upperK \upperA - 1] \sm A.
    \end{align*}
\end{Definition}

\begin{Proposition} \label{twoIneq}
	With the above notation, if for every $\upperMM$-avoiding set $\upperAA \subset \mathbb{Z}$ such that $0\in \upperAA$,
	\begin{align} \label{eq_twoIneq}
		\begin{aligned}
			\mbox{either } \  \upperA - |\upperII| + |\upperTT| &\le
			\begin{cases}
				 \upperA + \upperM \upperD, & 1 \le \upperR \le \upperK \\
				 \upperB - (\upperK+1)(\upperD+1), & \upperK+1 \le \upperR \le \upperM+\upperK,
			\end{cases} \\
			\mbox{or } \ \upperB - |\upperUU| - |\upperII| + |\upperTT| &\le
			\begin{cases} 
				\upperA + (\upperM+1) \upperD, & 1 \le \upperR \le \upperK \\
				\upperB - \upperK(\upperD+1), & \upperK+1 \le \upperR \le \upperM+\upperK,
			\end{cases}
		\end{aligned}
	\end{align}
	then the lower bound on $\mu(M)$ from the statement of \Cref{main-conjecture} is tight.
\end{Proposition}
\begin{proof}
	We need the following result from \cite[Lemma~1]{HARALAMBIS197722}, which is basically a simple averaging argument.
	
	\begin{Lemma}\label{haralambis}
		Let $\upperMM \subset \N$ be finite, and $\delta \in \mathbb{R}$ be positive. If for every $\upperMM$-avoiding set $\upperAA \subset \Z$ with $0\in \upperAA$, there exists $n \in \N$ such that $|\upperAA \cap [0,\upperN-1]| \le \delta\upperN$, then $\mu(\upperMM) \le \delta$.
	\end{Lemma}

	We claim that our set $M$ meets the conditions of \Cref{haralambis} if $\delta$ equals the lower bound on $\mu(M)$ from the statement of \Cref{main-conjecture} and $\upperN$ equals either $\upperN_1= \upperM \upperB + (\upperK + 1)\upperA$ or $\upperN_2=(\upperM + 1)\upperB + \upperK \upperA$.
	
	To confirm this claim, we first observe that
	\begin{equation} \label{eq_n1}
		[0,\upperN_1-1] = \bigg(\bigcup_{\upperAlpha = 0}^{\upperA-1} \upperSS_\upperAlpha \bigg) \cup \bigg(\bigcup_{i = 0}^{\upperM-1} [i\upperB + \upperA + 1, (i + 1)\upperB - 1] \bigg) \cup \{\upperA, \upperB+\upperA, \dots , (\upperM-1)\upperB+\upperA\}
	\end{equation} 
	and this union is disjoint. Moreover, the last set in \eqref{eq_n1} is disjoint from the $\upperMM$-avoiding set  $\upperAA$ since $0 \in \upperAA$. Besides that, each $\upperSS_\upperAlpha$ contains at most one element of $\upperAA$, and there are precisely $|\upperII|$ indices $\upperAlpha \in [0, \upperA-1]$ such that $\upperSS_\upperAlpha$ is disjoint from $\upperAA$ by the definition of $\upperII$. Hence, \eqref{eq_n1} implies that
	\begin{align*}
		|\upperAA \cap [0,\upperN_1-1]| &= \upperA-|\upperII|+|\upperTT|, \\
		|\upperAA \cap [0,\upperN_2-1]| &= |\upperAA \cap [0,\upperN_1-1]| + |\upperAA \cap [\upperN_1,\upperN_2-1]| \\ 
		&=  \upperA-|\upperII|+|\upperTT| + (\upperB-\upperA) - |\upperUU| \\
		&= \upperB-|\upperUU|-|\upperII|+|\upperTT|.
	\end{align*}
	We divide the remainder of this proof into two\footnote{As mentioned in the introduction, \Cref{main-conjecture} is trivially tight in case $r=0$.} cases depending on the value of $r$.
	
	First, suppose that $1 \le \upperR \le \upperK$. In this case, we have $\delta = \frac{\upperA + \upperM \upperD}{\upperN_1}$. If
	\begin{equation*}
		|\upperAA \cap [0,\upperN_1-1]| = \upperA-|\upperII|+|\upperTT| \le \upperA + \upperM \upperD = \delta\upperN_1,
	\end{equation*}
	then we are done. Otherwise, the assumption of \Cref{twoIneq} yields that
	\begin{equation*}
		|\upperAA \cap [0,\upperN_2-1]| = \upperB-|\upperUU|-|\upperII|+|\upperTT| \le \upperA + (\upperM+1) \upperD,
	\end{equation*}
	and it remains only to check that $\upperA + (\upperM+1) \upperD < \delta\upperN_2$, because
	\begin{equation*}
		(\upperA + \upperM \upperD)\upperN_2 - (\upperA + (\upperM+1) \upperD)\upperN_1 = \upperA\upperR > 0.
	\end{equation*}

	Similarly, suppose that $\upperK+ 1 \le \upperR \le \upperM+\upperK$. In this case, we have $\delta = \frac{\upperB - \upperK(\upperD+1)}{\upperN_2}$. If
	\begin{equation*}
		|\upperAA \cap [0,\upperN_2-1]| = \upperB-|\upperUU|-|\upperII|+|\upperTT| \le \upperB - \upperK(\upperD+1) = \delta\upperN_2,
	\end{equation*}
	then we are done. Otherwise, the assumption of \Cref{twoIneq} yields that
	\begin{equation*}
		|\upperAA \cap [0,\upperN_1-1]| = \upperA-|\upperII|+|\upperTT| \le \upperB - (\upperK+1)(\upperD+1),
	\end{equation*}
	and it remains only to check that $\upperB - (\upperK+1)(\upperD+1) < \delta\upperN_1$, because
	\begin{equation*}
		(\upperB - \upperK(\upperD+1))\upperN_1 - (\upperB - (\upperK+1)(\upperD+1))\upperN_2 = \upperB(\upperM+\upperK+1-\upperR) > 0. \qedhere
	\end{equation*}
\end{proof}

\begin{Corollary} \label{oneIneq}
	Under the notation of Definition \ref{def_sets}, if
	\begin{equation}\label{weak-assumption}
		\upperM + \upperK >
		\begin{cases}
			\upperM\upperR, & 1 \le \upperR \le \upperK \\
			\upperK(\upperM+\upperK+1-\upperR), & \upperK+1 \le \upperR \le \upperM+\upperK,
		\end{cases}
	\end{equation}
	and for every $\upperMM$-avoiding set $\upperAA \subset \mathbb{Z}$ such that $0\in \upperAA$, we have
	\begin{equation}\label{main-inequality}
		(\upperM + \upperK) |\upperTT| \le (\upperM + \upperK) |\upperII| + \upperM|\upperUU|,
	\end{equation}
	then the lower bound on $\mu(M)$ from the statement of \Cref{main-conjecture} is tight.
\end{Corollary}
\begin{proof}
	It is sufficient to verify that \eqref{weak-assumption} and \eqref{main-inequality} imply \eqref{eq_twoIneq}, and then the conclusion follows from \Cref{twoIneq}. Therefore, let us assume the contrary.
	
	For $1 \le \upperR \le \upperK$, this assumption means that there exists an $\upperMM$-avoiding set $\upperAA \subset \mathbb{Z}$ such that $0\in \upperAA$ and
	\begin{equation*}
		\begin{cases*}
			\upperA - |\upperII| + |\upperTT| \ge \upperA + \upperM \upperD + 1, \\
			\upperB - |\upperUU| - |\upperII| + |\upperTT| \ge \upperA + (\upperM+1) \upperD + 1.
		\end{cases*}
	\end{equation*}
	We multiply the first inequality by $\upperK$, the second one by $\upperM$, and sum these products. After a straightforward rearrangement of terms, this leads to
	\begin{equation} \label{eq_cor1}
		(\upperM + \upperK) |\upperTT| - (\upperM + \upperK) |\upperII| -  \upperM|\upperUU| \ge \upperM + \upperK - \upperM\upperR.
	\end{equation}
	However, the right-hand side of \eqref{eq_cor1} is positive by \eqref{weak-assumption}, while the left-hand side is not by \eqref{main-inequality}, a contradiction.
	
	Similarly, for $\upperK+1 \le \upperR \le \upperM+\upperK$, our assumption means that there exists an $\upperMM$-avoiding set $\upperAA \subset \mathbb{Z}$ such that $0\in \upperAA$ and
	\begin{equation*}
		\begin{cases*}
			\upperA - |\upperII| + |\upperTT| \ge \upperB - (\upperK+1)(\upperD+1) + 1, \\
			\upperB - |\upperUU| - |\upperII| + |\upperTT| \ge \upperB - \upperK(\upperD+1) + 1.
		\end{cases*}
	\end{equation*}
	As before, we multiply the first inequality by $\upperK$, the second one by $\upperM$, and sum these products. After a straightforward rearrangement of terms, this leads to
	\begin{equation} \label{eq_cor2}
		(\upperM + \upperK) |\upperTT| - (\upperM + \upperK) |\upperII| -  \upperM|\upperUU| \ge \upperM + \upperK - \upperK(\upperM+\upperK+1-\upperR).
	\end{equation}
	Again, the right-hand side of \eqref{eq_cor2} is positive by \eqref{weak-assumption}, while the left-hand side is not by \eqref{main-inequality}, a contradiction.
\end{proof}

Note that \eqref{weak-assumption} holds for all $\upperR$ whenever either $\upperM$ or $\upperK$ equals $1$. In the next two subsections, we show that \eqref{main-inequality} holds in these two special cases as well, so the conditions of \Cref{oneIneq} are met.

\subsection[The first case]{Auxiliary inequality --- proof for $\upperK = 1$}

In this subsection, we suppose that $\upperK = 1$ and fix an arbitrary $M$-avoiding set $\upperAA \subset \Z$ such that $0 \in \upperAA$. Our goal is to show that the inequality \eqref{main-inequality} holds in this case, i.e., that
\begin{equation}\label{main-inequality-m=1}
	(\upperM + 1) |\upperTT| \le (\upperM + 1) |\upperII| + \upperM|\upperUU|,
\end{equation}
and then \Cref{oneIneq} will complete the proof of \Cref{main-theorem} for the case $\upperB > \upperA$.

To this end, we construct a partition $\upperTT = \bigcup_{\upperI} \upperCC_{\upperI}$ of $\upperTT$ into the union of disjoint `chains' and map each chain into a subset of $\upperII \cup \upperUU$ such that the following three properties hold:
\begin{property}
    \item \label{chain-part-not-exceed-k} the size of each $\upperCC_{\upperI}$ satisfies $|\upperCC_{\upperI}| \le \upperM$,
    \item\label{number-of-corresponded-points} the image of each $\upperCC_{\upperI}$ is of size $|\upperCC_{\upperI}| + 1$,
    \item\label{injectivity} the images of different chains are disjoint.
\end{property}
Since the function $\frac{x}{x+1}$ is increasing on $\N$, the first two properties imply that the size of each $\upperCC_i \subset \upperTT$ is at most $\frac{\upperM}{\upperM+1}$ times the size of its image in $\upperII \cup \upperUU$. Together with the third property, this yields that
\begin{equation*}
	|\upperTT| \le \frac{\upperM}{\upperM+1}|\upperII\cup\upperUU|,
\end{equation*}
which is stronger than the desired inequality \eqref{main-inequality-m=1}. To construct this partition, we first map each $\upperAlpha \in T$ to a  two-element set $\{\upperV_\upperAlpha, \upperW_\upperAlpha\} \subset \upperII \cup \upperUU$, and then we analyze possible intersections of these sets.

%let $\upperQ, \upperP, \upperEta, \widehat \upperAlpha_0$ be unique integers such that
%\begin{align*}
    %& \upperB = \upperQ\upperA + \upperP, \qquad \quad \;\; 0 \le \upperP < \upperA,\\
    %& \upperAlpha = \upperXi \upperB + \upperEta \upperA + \widehat \upperAlpha_0, \;\; 1 \le \upperEta \le \upperQ, \;\; 0 \le \widehat \upperAlpha_0 < \upperA.
%\end{align*}
Given $0 \le \upperXi < \upperM$, $\upperAlpha \in \upperTT_\upperXi = \upperAA \cap [\upperXi\upperB + \upperA + 1, (\upperXi + 1)\upperB - 1]$, let $\upperEta \ge 1$, $0 \le \widehat \upperAlpha_0 < \upperA$ be the unique integers such that $\upperAlpha = \upperXi \upperB + \upperEta \upperA + \widehat \upperAlpha_0$. Define $\upperV_\upperAlpha$ as $\upperAlpha + (\upperM - \upperXi)\upperB + \upperA$ and note that $\upperV_\upperAlpha-\upperAlpha \in \upperMM$, so $\upperV_\upperAlpha \notin A$. Besides, we have
\begin{equation*}
	\upperM \upperB + 2\upperA \le \upperM \upperB + (\upperEta+1) \upperA + \widehat \upperAlpha_0 = \upperV_\upperAlpha = \upperAlpha + (\upperM - \upperXi)\upperB + \upperA \le  (\upperM + 1)\upperB + \upperA - 1,
\end{equation*}
and thus $\upperV_\upperAlpha \in \upperUU = [\upperM \upperB + 2\upperA, (\upperM + 1)\upperB + \upperA - 1] \sm A.$

Similarly, if $\upperEta \ge 2$, we put $\upperW_\upperAlpha= \upperAlpha + (\upperM - \upperXi)\upperB$ and note that $\upperW_\upperAlpha \in \upperUU$ since
\begin{equation*}
	\upperM \upperB + 2\upperA \le \upperM \upperB + \upperEta \upperA + \widehat \upperAlpha_0 = \upperW_\upperAlpha = \upperAlpha + (\upperM - \upperXi)\upperB \le  (\upperM + 1)\upperB - 1.
\end{equation*}
%Observe that if $\upperP=0$, then $\upperEta \ge 2$, and we have already defined $\upperW_\upperAlpha$. Indeed, in this case, we have $\upperB = \upperQ\upperA$, and thus $\upperA=1$ since $\upperB$ and $\upperA$ are coprime. If, in addition, we have $\upperEta = 1$, then $\upperAlpha-0 = \upperXi \upperB + \upperA \in \upperMM$, a contradiction. Therefore, it remains only to define $\upperW_\upperAlpha$ under the additional assumptions $\upperP>0$ and $\upperEta = 1$.

In the remaining case $\upperEta=1$, we consider a sequence $\widehat\upperAlpha_\upperI = \widehat\upperAlpha_0+\upperI(\upperB-\upperA)$, $\upperI \in \N$. Let $\upperNN$ be the minimum integer such that either $\widehat\upperAlpha_\upperNN \in \upperII$ or $\widehat\upperAlpha_\upperNN \ge 2\upperA-\upperB$. We define $\upperW_\upperAlpha$ as $\widehat\upperAlpha_\upperNN$ if $\widehat\upperAlpha_\upperNN \in \upperII$ and as $\widehat\upperAlpha_\upperNN + (\upperM + 1)\upperB$ if $\widehat\upperAlpha_\upperNN \notin \upperII$. Our proof that $\upperW_\upperAlpha \in \upperUU$ for the case $\widehat \upperAlpha_\upperNN \ge 2\upperA-\upperB$ relies on the following result.

\begin{Lemma}\label{late-sequence}
	For all $\upperAlpha = \upperXi \upperB + \upperA + \widehat \upperAlpha_0 \in \upperTT_{\upperXi}$ and $0 \le \upperI \le \upperNN$, if $\widehat \upperAlpha_\upperI \notin \upperII$, then $\widehat \upperAlpha_\upperI + \upperM' \upperB \in \upperAA$ for some $\upperXi < \upperM' \le \upperM$.
\end{Lemma}
\begin{proof}
	We argue by induction on $\upperI$. Note that the condition $\widehat \upperAlpha_0 \not\in \upperII$ yields that $\upperAA \cap \upperSS_{\widehat \upperAlpha_0} \neq \varnothing$, where  $\upperSS_{\widehat \upperAlpha_0} = \{\widehat \upperAlpha_0, \widehat \upperAlpha_0 + \upperB,  \dots, \widehat \upperAlpha_0 + \upperM\upperB, \widehat \upperAlpha_0 + \upperM\upperB + \upperA\}$. Hence, it is sufficient to prove that the first $\upperXi + 1$ elements of $\upperSS_{\widehat \upperAlpha_0}$ and the last one are not in $\upperAA$. If $\widehat \upperAlpha_0 + \upperM \upperB + \upperA \in \upperAA$, we have $(\widehat \upperAlpha_0 + \upperM \upperB + \upperA) - \upperAlpha = (\upperM - \upperXi) \upperB \in \upperMM$, a contradiction. If $\widehat \upperAlpha_0 + \upperM'\upperB \in \upperAA$, $0 \le \upperM' \le \upperXi$, then $\upperAlpha - (\widehat \upperAlpha_0 + \upperM'\upperB) = (\upperXi - \upperM')\upperB + \upperA \in \upperMM$, a contradiction again.
	
	Now we prove the induction step. Suppose that $\upperI<\upperNN$, and take $\upperXi < \upperM' \le \upperM$ such that $\widehat \upperAlpha_\upperI + \upperM' \upperB \in \upperAA$, which exists by the induction hypothesis. As before, the condition $\widehat \upperAlpha_{\upperI+1} \not\in \upperII$ yields that $\upperAA \cap \upperSS_{\widehat \upperAlpha_{\upperI+1}} \neq \varnothing$, so our goal is to prove that neither the first $\upperXi + 1$ elements of $\upperSS_{\widehat \upperAlpha_{\upperI+1}}$ are in $\upperAA$ nor the last one is. If $\widehat \upperAlpha_{\upperI+1} + \upperM \upperB + \upperA \in \upperAA$, then $(\widehat \upperAlpha_{\upperI+1} + \upperM \upperB + \upperA) - (\widehat \upperAlpha_{\upperI} + \upperM'\upperB) = (\upperM - \upperM' + 1)\upperB \in \upperMM$, a contradiction. If $\widehat \upperAlpha_{\upperI+1} + \upperM'' \upperB \in \upperAA$ for some $0 \le \upperM'' \le \upperXi$, then $(\widehat \upperAlpha_{\upperI} + \upperM'\upperB) - (\widehat \upperAlpha_{\upperI+1} + \upperM''\upperB) = (\upperM' - \upperM'' - 1)\upperB + \upperA \in \upperMM$, a contradiction again.
\end{proof}

Now let us check that if $\upperAlpha = \upperXi \upperB + \upperA + \widehat \upperAlpha_0 \in \upperTT_\upperXi$ and $\upperW_\upperAlpha = \widehat\upperAlpha_\upperNN + (\upperM + 1)\upperB$, then $\upperW_\upperAlpha \in \upperUU$. First, note that
\begin{equation*}
	\upperM \upperB + 2\upperA = (2\upperA-\upperB) + (\upperM + 1)\upperB  \le \widehat\upperAlpha_\upperNN + (\upperM + 1)\upperB \le  (\upperM + 1)\upperB + \upperA - 1.
\end{equation*}
Moreover, we have $\widehat\upperAlpha_\upperNN \notin \upperII$ by construction, and thus $\widehat \upperAlpha_\upperNN + \upperM'\upperB \in \upperAA$ for some $\upperXi < \upperM' \le \upperM$ by \Cref{late-sequence}. Since $\upperW_\upperAlpha - (\widehat \upperAlpha_\upperNN + \upperM'\upperB) = (\upperM-\upperM' + 1)\upperB \in \upperMM$, we conclude that $\upperW_\upperAlpha \notin \upperAA$ and thus $\upperW_\upperAlpha \in \upperUU$, as desired. Besides, observe that $\upperW_\upperAlpha \neq \upperV_\upperAlpha$ since otherwise $(\widehat \upperAlpha_\upperNN + \upperM' \upperB) - \upperAlpha = (\upperM' - \upperXi - 1)\upperB + \upperA \in \upperMM$, a contradiction.

Therefore, we indeed map each $\upperAlpha \in \upperTT$ to a well-defined two-element subset $\{\upperV_\upperAlpha, \upperW_\upperAlpha\} \subset \upperII \cup \upperUU$. Next, let us study potential intersections between these subsets.

%\begin{Lemma}\label{alpha-sequence}
%	If $\upperAlpha = \upperXi \upperB + \upperA + \widehat \upperAlpha_0$ and $\upperW_\upperAlpha = \widehat\upperAlpha_\upperNN + (\upperM + 1)\upperB$, then $\upperW_\upperAlpha \in \upperUU$.
%\end{Lemma}

%\begin{proof}
	%The translate $\upperSS_{\widehat \upperAlpha_{0}}$ is $\upperXi$-late by \Cref{explicit-intersecting-constant}. Note that\footnote{The corresponding argument from \cite[Section~4]{liuzhu2004} was valid only in this special case, and was later extended for all $\upperB, \upperA$ in \cite{liu2005d}.} if $\upperB\ge 2\upperA$, then $\upperNN=0$. Otherwise, we have $\upperQ=1$ and $\widehat\upperAlpha_\upperI = \widehat\upperAlpha_{i-1}+\upperB-\upperA$ for all $\upperI \in [\upperNN]$. Now \Cref{invariant} applied $\upperNN$ times yields that $\upperSS_{\widehat\upperAlpha_\upperNN}$ is $\upperXi$-late. In either case, we have $\widehat \upperAlpha_\upperNN + \upperM'\upperB \in \upperAA$, for some $\upperXi < \upperM' \le \upperM$. Since $\upperW_\upperAlpha - (\widehat \upperAlpha_\upperNN + \upperM'\upperB) = (\upperM-\upperM' + 1)\upperB \in \upperMM$, we conclude that $\upperW_\upperAlpha \notin \upperAA$. Moreover, it is not hard to check that
	%\begin{equation*}
	%	\upperM \upperB + 2\upperA = (2\upperA-\upperB) + (\upperM + 1)\upperB  \le \widehat\upperAlpha_\upperNN + (\upperM + 1)\upperB \le  (\upperM + 1)\upperB + \upperA - 1,
	%\end{equation*}
	%and thus $\upperW_\upperAlpha = \widehat \upperAlpha_\upperNN + (\upperM + 1)\upperB \in \upperUU$, as desired.
%\end{proof}

\begin{Lemma}\label{t-chain-neighbor-criterion}
	\leavevmode
	\begin{itemize}
		\item For each $\upperAlpha \in \upperTT_{\upperXi}$, there exists at most one $\upperBeta \in \upperTT$ such that $\upperAlpha > \upperBeta$ and $\{\upperV_\upperAlpha, \upperW_\upperAlpha\} \cap \{\upperV_\upperBeta, \upperW_\upperBeta\} \neq \varnothing$. Moreover, for such $\upperBeta$, this intersection is of size $1$ and $\upperBeta \notin \upperTT_{\upperXi}$.
		\item For each $\upperBeta \in \upperTT_{\upperXi'}$, there exists at most one $\upperAlpha \in \upperTT$ such that $\upperAlpha > \upperBeta$ and $\{\upperV_\upperAlpha, \upperW_\upperAlpha\} \cap \{\upperV_\upperBeta, \upperW_\upperBeta\} \neq \varnothing$. Moreover, for such $\upperAlpha$, this intersection is of size $1$ and $\upperAlpha \notin \upperTT_{\upperXi'}$.
	\end{itemize}
\end{Lemma}

\begin{proof}
	Fix $\upperAlpha \in \upperTT_{\upperXi}$ and $\upperBeta \in \upperTT_{\upperXi'}$ for some $0 \le \upperXi' \le \upperXi < \upperM$. By construction, we have  $\upperV_\upperAlpha = \upperAlpha + (\upperM - \upperXi)\upperB+ \upperA$ and  $\upperW_\upperAlpha \in \{\upperAlpha + (\upperM - \upperXi)\upperB, \widehat \upperAlpha_N + (\upperM + 1)\upperB, \widehat \upperAlpha_N\}$. Similarly, $\upperV_\upperBeta = \upperBeta + (\upperM - \upperXi')\upperB+ \upperA$ and  $\upperW_\upperBeta \in \{\upperBeta + (\upperM - \upperXi')\upperB, \widehat \upperBeta_{\upperNN'} + (\upperM + 1)\upperB, \widehat \upperBeta_{\upperNN'}\}$. Let us separately consider all potential ways for some two of these numbers to coincide.
	\begin{enumerate}
		\item If $\upperAlpha \neq \upperBeta$ and either $\upperAlpha + (\upperM - \upperXi)\upperB + \upperA = \upperBeta + (\upperM - \upperXi')\upperB + \upperA$ or $\upperAlpha + (\upperM - \upperXi) \upperB = \upperBeta + (\upperM - \upperXi') \upperB$, then  $\upperAlpha - \upperBeta = (\upperXi - \upperXi')\upperB \in \upperMM$, a contradiction.
		\item If $\upperAlpha + (\upperM - \upperXi) \upperB = \upperBeta + (\upperM - \upperXi')\upperB + \upperA$, then $\upperAlpha-\upperBeta =  (\upperXi - \upperXi')\upperB + \upperA \in \upperMM$, a contradiction.
		\item If $\upperAlpha + (\upperM - \upperXi) \upperB = \widehat \upperBeta_{\upperNN'} + (\upperM + 1)\upperB$, we get $\upperAlpha = \widehat\upperBeta_{\upperNN'} + (\upperXi + 1)\upperB  \ge (\upperXi + 1)\upperB$ and thus $\upperAlpha \notin \upperTT_\upperXi$, a contradiction.
		\item Similarly, if $\widehat\upperAlpha_\upperNN + (\upperM + 1)\upperB = \upperBeta + (\upperM - \upperXi')\upperB$, then $\upperBeta = \widehat\upperAlpha_\upperNN + (\upperXi' + 1)\upperB \not\in \upperTT_{\upperXi'}$, a contradiction.
		\item If $\upperW_\upperAlpha = \widehat\upperAlpha_\upperNN + (\upperM + 1)\upperB = \upperBeta + (\upperM - \upperXi')\upperB + \upperA$, \Cref{late-sequence} implies that $\widehat\upperAlpha_\upperNN + \upperM' \upperB \in \upperAA$ for some $\upperXi < \upperM' \le \upperM$, and thus $(\widehat \upperAlpha_\upperNN + \upperM' \upperB) - \upperBeta = (\upperM' - \upperXi' - 1)\upperB + \upperA \in \upperMM$, a contradiction.
		\item If either $\upperW_\upperAlpha = \widehat\upperAlpha_\upperNN + (\upperM + 1)\upperB = \widehat\upperBeta_{\upperNN'} + (\upperM + 1)\upperB = \upperW_\upperBeta $ or $\upperW_\upperAlpha = \widehat\upperAlpha_\upperNN = \widehat\upperBeta_{\upperNN'} = \upperW_\upperBeta $, then $\upperAlpha = \upperXi \upperB + \upperA + \widehat \upperAlpha_0$ and $\upperBeta = \upperXi' \upperB + \upperA + \widehat \upperBeta_0$ by construction. Without loss of generality, assume that $\upperNN \ge \upperNN'$. Then the equality $\widehat\upperAlpha_\upperNN = \widehat\upperBeta_{\upperNN'}$ implies that $\widehat\upperAlpha_\upperI = \widehat\upperBeta_{0}$, where $\upperI=\upperNN - \upperNN'$. If $\upperI=0$, then $\upperAlpha - \upperBeta = (\upperXi - \upperXi')\upperB \in \upperMM$ whenever $\upperAlpha \neq \upperBeta$, a contradiction. Otherwise, we have and $\widehat\upperAlpha_\upperI = \widehat\upperAlpha_{\upperI-1}+\upperB-\upperA$ by construction. Hence, $\upperBeta = (\upperXi'+1)\upperB+\widehat\upperAlpha_{\upperI-1} \ge (\upperXi'+1)\upperB$ and thus $\upperBeta \notin \upperTT_{\upperXi'}$, a contradiction again.
	\end{enumerate}
	The only coincidence that we have not excluded yet is $\upperV_\upperAlpha=\upperW_\upperBeta$, and it indeed might occur. However, observe that for each fixed $\upperBeta \in \upperTT$, the equality $\upperV_\upperAlpha=\upperW_\upperBeta$ can hold for at most one $\upperAlpha \in \upperTT$ such that $\upperAlpha > \upperBeta$. Indeed, if $\upperV_\upperAlpha=\upperW_\upperBeta=\upperV_{\upperAlpha'}$ for some $\upperAlpha > \upperAlpha' > \upperBeta$, then the first out of six statements above applied to $\upperAlpha'$ playing the role of $\upperBeta$ yields a contradiction.  Similarly, if $\upperW_\upperBeta=\upperV_\upperAlpha=\upperW_{\upperBeta'}$ for some $\upperAlpha > \upperBeta' > \upperBeta$, then one of the six statements above applied to $\upperBeta'$ playing the role of $\upperAlpha$ yields a contradiction. 
	
	Hence, to complete the proof, it remains only to show that if $\upperV_\upperAlpha=\upperW_\upperBeta$, then $\upperXi>\upperXi'$. We consider the following two cases. First, suppose that $\upperW_\upperBeta = \upperBeta + (\upperM - \upperXi')\upperB$. Then $\upperV_\upperAlpha=\upperW_\upperBeta$ implies that
	$0 \le \upperAlpha - \upperBeta =  (\upperXi - \upperXi')\upperB - \upperA$, and thus $\upperXi>\upperXi'$, as desired. Second, suppose that $ \upperW_\upperBeta = \widehat\upperBeta_{\upperNN'} + (\upperM + 1) \upperB$. Then $\upperBeta = \upperXi' \upperB + \upperA + \widehat\upperBeta_0$ and $\widehat\upperBeta_{\upperNN'} + \upperM' \upperB \in \upperAA$ for some $\upperXi' < \upperM' \le \upperM$ by \Cref{late-sequence}. If $\upperXi=\upperXi'$, then $\upperM' - \upperXi - 1 \ge 0$. Now $\upperV_\upperAlpha=\upperW_\upperBeta$ implies that $(\widehat\upperBeta_{\upperNN'} + \upperM' \upperB)- \upperAlpha = (\upperM' - \upperXi - 1)\upperB + \upperA \in \upperMM$, a contradiction.
\end{proof}

Now we use this result to construct the desired partition $\upperTT =  \bigcup_{\upperI} \upperCC_{\upperI}$. Consider a directed graph on the vertex set $\upperTT$, where $(\upperAlpha, \upperBeta)$ is an edge if and only if $\upperAlpha > \upperBeta$ and $\{\upperV_\upperAlpha, \upperW_\upperAlpha\} \cap \{\upperV_\upperBeta, \upperW_\upperBeta\} \neq \varnothing$. \Cref{t-chain-neighbor-criterion} implies that both the indegree and the outdegree of each vertex are at most $1$. Hence, our directed graph is a disjoint union of ordered path that we call `chains'. Naturally, we map each chain $\upperCC_{\upperI}$ into $\bigcup_{\upperAlpha \in \upperCC_{\upperI}} \{\upperV_\upperAlpha, \upperW_\upperAlpha\} \subset \upperII\cup\upperUU$. Since there are no edges between different chains, their images are disjoint, and so the property \eqref{injectivity} holds. Besides that, for each chain, the images of two distinct vertices on it either have precisely one common point if these vertices are consecutive or are  otherwise disjoint by \Cref{t-chain-neighbor-criterion}. Now it is not hard to see that the size of each chain is less than the size of its image by exactly $1$, and so \eqref{number-of-corresponded-points} holds as well. Finally, observe that the index $0 \le \upperXi < \upperM$ corresponding to each vertex $\upperAlpha \in \upperTT_{\upperXi}$ strictly decreases as we go along the  edges of a chain by \Cref{t-chain-neighbor-criterion}. Therefore, the size of chain does not exceed $\upperM$, the property \eqref{chain-part-not-exceed-k} holds, and we are done.

\subsection[The second case]{Auxiliary inequality --- proof for $\upperM = 1$} \label{Sec3.3}

In this subsection, we prove \Cref{main-theorem} in the remaining case $\upperB < \upperA$. First, we utilize the symmetric nature of this problem discussed in the introduction to make the results from \Cref{Sec3.1} applicable in this case: replace $\upperSS$ with its reflection and switch the roles of $\upperB$ and $\upperA$, as well as of $\upperM$ and $\upperK$. In other words, we suppose that $\upperB > \upperA$ and
\begin{equation*}
	S = \{0, \theora, \theora + \theorb, \dots, \theora + \theorm\theorb\}, \ \ \ 
	M=\{\upperA, \dots, \upperK\upperA\}\cup \{\upperB, \upperB+\upperA, \dots, \upperB+\upperK\upperA\}.
\end{equation*}
Now a tight upper bound on $\mu(M)$ follows from \Cref{oneIneq} if the inequality \eqref{main-inequality} holds for $\upperM = 1$, i.e., if
\begin{equation*}\label{main-inequality-k=1}
	(\upperK + 1) |\upperTT| \le (\upperK + 1) |\upperII| + |\upperUU|,
\end{equation*}
for every $M$-avoiding set $\upperAA \subset \Z$ such that $0 \in \upperAA$. To establish the inequality above, we fix one such set $\upperAA$ and map each point $\upperAlpha \in \upperTT$ either into one element in $\upperII$ or into $(\upperK + 1)$-element subset $\upperUU(\upperAlpha) \subset \upperUU$ in such a way that the images of distinct points are disjoint. This mapping shares many similarities with the one constructed in the previous subsection.

%As earlier, given $\upperAlpha \in \upperTT = \upperAA \cap [\upperA + 1, \upperB - 1]$, let $\upperQ, \upperP, \upperEta_0, \widehat \upperAlpha_0$ be unique integers such that
%\begin{align*}
	%& \upperB = \upperQ\upperA + \upperP, \quad \; 0 \le \upperP < \upperA,\\
	%& \upperAlpha = \upperEta_0 \upperA + \widehat \upperAlpha_0, \;\; 1 \le \upperEta_0 \le \upperQ, \;\; 0 \le \widehat \upperAlpha_0 < \upperA.
%\end{align*}

As before, given $\upperAlpha \in \upperTT = \upperAA \cap [\upperA + 1, \upperB - 1]$, let $\upperEta_0\ge 1, 0 \le \widehat \upperAlpha_0 < \upperA$ be the unique integers such that $\upperAlpha = \upperEta_0 \upperA + \widehat \upperAlpha_0$.

If $\upperEta_0 >  \upperK$, we define 
\begin{equation*}
	\upperUU(\upperAlpha) = \upperVV (\upperAlpha) = \{\upperAlpha+\upperB +\upperJ\upperA: 0 \le \upperJ \le \upperK \}.
\end{equation*}
%Observe that if $\upperP=0$, then $\upperEta_0 > \upperK$. Indeed, in this case, we have $\upperB = \upperQ\upperA$, and thus $\upperA=1$ since $\upperB$ and $\upperA$ are coprime. If, in addition, we have $\upperEta_0 \le \upperK$, then $\upperAlpha-0 = \upperEta_0\upperA \in \upperMM$, a contradiction. Therefore, it remains only to define $\upperUU(\upperAlpha)$ under the additional assumptions $\upperP>0$ and $\upperEta_0 \le \upperK$. 

%We make a few observations that are immediate from this definition:
%\begin{observation}
%\item \label{obs-N+1vsN} $ \upperTau_{\upperN+1}\upperA+ \widehat \upperAlpha_{\upperN + 1}=  \upperB+(\upperTau_{\upperN}-1)\upperA+\widehat \upperAlpha_{\upperN}$,
%\end{observation}

Assume now that $\upperEta_0 \le \upperK$. For $\upperN \in \N$, define  $\upperTau_{\upperN}, \widehat \upperAlpha_\upperN$ as the unique integers such that
\begin{equation*}
	\upperAlpha+\upperN(\upperB-\upperA)=\upperTau_{\upperN} \upperA + \widehat \upperAlpha_\upperN, \;\; 0 \le \widehat \upperAlpha_{\upperN} < \upperA
\end{equation*}
and observe that
\begin{equation} \label{obs-N+1vsN}
	\upperTau_{\upperN+1}\upperA+ \widehat \upperAlpha_{\upperN + 1}=  \upperB+(\upperTau_{\upperN}-1)\upperA+\widehat \upperAlpha_{\upperN}.
\end{equation}
Let $\upperNN$ be the minimum integer such that either $\widehat\upperAlpha_\upperNN \in \upperII$ or $\upperTau_{\upperNN+1} \ge \upperK+1$. If $\widehat\upperAlpha_\upperNN \in \upperII$, then we map $\upperAlpha$ into $\widehat\upperAlpha_\upperNN$. Otherwise $\upperTau_{\upperNN+1} \ge \upperK+1$, and we `truncate' $\upperTau_{\upperNN+1}$ by setting $\upperTau_{\upperNN+1} = \upperK+1$. Next, we define
\begin{align*} \label{alpha-sequence-v}
	\begin{aligned}
		\upperVV (\upperAlpha) &= \{\upperAlpha + \upperB + \upperJ\upperA : \upperK - \upperEta_0 + 1 \le \upperJ \le \upperK\},\\
		\upperWW_\upperN (\upperAlpha) &= \{\widehat \upperAlpha_\upperN + 2\upperB + \upperJ\upperA : \upperK + \upperTau_{\upperN} -\upperTau_{\upperN+1} \le \upperJ < \upperK \}, \;\; 0 \le \upperN \le \upperNN, \phantom{\bigcup^{\upperNN}} \\
		\upperUU(\upperAlpha) &= \upperVV (\upperAlpha) \cup \bigcup_{\upperN = 0}^{\upperNN} \upperWW_\upperN(\upperAlpha),
	\end{aligned}
\end{align*}
and map $\upperAlpha$ into $\upperUU (\upperAlpha)$. Our proof that $\upperUU (\upperAlpha) \subset \upperUU$ relies on the following result.

\begin{Lemma}\label{t-for-main-sequence-lemma}
	For all $\upperAlpha  = \upperEta_0 \upperA + \widehat \upperAlpha_0 \in \upperTT$ such that $1 \le \upperEta_0 \le \upperK$ and for all $0 \le \upperN \le \upperNN$, if $\widehat \upperAlpha_\upperN \not\in \upperII$, then $\widehat \upperAlpha_{\upperN} + \upperB + \upperK'\upperA \in \upperAA$ for some $0 \le \upperK' < \upperTau_{\upperN}$.
\end{Lemma}
\begin{proof}
	We argue by induction on $\upperN$. Note that the condition $\widehat \upperAlpha_0 \not\in \upperII$ yields that  $\upperAA \cap \upperSS_{\widehat \upperAlpha_0} \neq \varnothing$, where $\upperSS_{\widehat \upperAlpha_0} = \{\widehat \upperAlpha_0, \widehat \upperAlpha_0 + \upperB, \widehat \upperAlpha_0 + \upperB + \upperA, \dots, \widehat \upperAlpha_0 + \upperB + \upperK \upperA\}$. Hence, it is sufficient to prove that the first one and the last $\upperK - \upperEta_0 + 1$ elements of $S_{\widehat\upperAlpha_0}$ are not in $\upperAA$. For the first one, we have $\upperAlpha - \widehat\upperAlpha_0 = \upperEta_0 \upperA \in \upperMM$, so $\widehat\upperAlpha_0 \not\in \upperAA$. Similarly, if $\upperEta_0 \le \upperK' \le \upperK$, then $(\widehat\upperAlpha_0 + \upperB + \upperK'\upperA) - \upperAlpha = \upperB + (\upperK' - \upperEta_0)\upperA \in \upperMM$, so $\widehat\upperAlpha_0 + \upperB + \upperK'\upperA \notin \upperAA$, as desired.
	
	Now we prove the induction step. Suppose the lemma is true for $\upperN < \upperNN$, i.e., that $\widehat\upperAlpha_\upperN + \upperB + \upperK'\upperA \in \upperAA$ for some $0 \le \upperK' < \upperTau_{\upperN}$. Note that $\upperN< \upperNN$ implies $\upperTau_{\upperN + 1} \le \upperK$. As earlier, the condition $\widehat \upperAlpha_{\upperN+1} \not\in \upperII$ yields that $\upperAA \cap \upperSS_{\widehat \upperAlpha_{\upperN+1}} \neq \varnothing$, so our goal is to prove that neither the first one nor the last $\upperK - \upperEta_{\upperN+1} + 1$ elements of $S_{\widehat\upperAlpha_{\upperN+1}}$ are in $\upperAA$. First, we have $1 \le (\upperTau_{\upperN + 1}-\upperTau_{\upperN})+\upperK'+1 = \upperTau_{\upperN + 1}-(\upperTau_{\upperN}-\upperK'-1) \le \upperK.$ Therefore, \eqref{obs-N+1vsN} implies  that $(\widehat\upperAlpha_\upperN + \upperB + \upperK'\upperA)-\widehat\upperAlpha_{\upperN+1} = (\upperTau_{\upperN + 1}-\upperTau_{\upperN}+\upperK'+1)\upperA \in \upperMM$, so $\widehat\upperAlpha_{\upperN+1} \notin \upperAA$. Similarly, if $\upperEta_{\upperN+1} \le \upperK'' \le \upperK$, then $0 \le (\upperK''-\upperTau_{\upperN + 1})+(\upperTau_{\upperN}-\upperK'-1) = \upperK''-(\upperTau_{\upperN + 1}-\upperTau_{\upperN})-(\upperK'+1) \le \upperK$. Hence, it follows from \eqref{obs-N+1vsN} that $(\widehat\upperAlpha_{\upperN+1} + \upperB + \upperK''\upperA)-(\widehat\upperAlpha_\upperN + \upperB + \upperK'\upperA) = \upperB+(\upperK''-\upperTau_{\upperN + 1}+\upperTau_{\upperN}-\upperK'-1)\upperA \in \upperMM$, and thus $\widehat\upperAlpha_{\upperN+1} + \upperB + \upperK''\upperA \notin \upperAA$.
\end{proof}

\begin{Lemma}\label{u-alpha-is-in-u}
	If $\upperAlpha\in \upperTT$ is not mapped into $\widehat \upperAlpha_\upperNN$, then $\upperUU (\upperAlpha) \subset \upperUU$.
\end{Lemma}

\begin{proof}
	Fix some $\upperAlpha  = \upperEta_0 \upperA + \widehat \upperAlpha_0 \in \upperTT$, and recall that $\upperUU = [\upperB + (\upperK + 1)\upperA, 2\upperB + \upperK \upperA-1]\sm \upperAA$.
	
	First, we prove that $\upperVV (\upperAlpha) \subset \upperUU$. Take an arbitrary element $\upperAlpha + \upperB + \upperJ\upperA$ of $\upperVV (\upperAlpha)$. Since $0 \le \upperJ \le \upperK$ and $\upperJ \ge \upperK - \upperEta_0 + 1$, we have $(\upperAlpha + \upperB + \upperJ\upperA)-\upperAlpha \in \upperMM$, so $\upperAlpha + \upperB + \upperJ\upperA \notin \upperAA$. Besides, we have
	\begin{equation*}
		\upperB + (\upperK + 1)\upperA \le \upperB + (\upperEta_0 + \upperJ)\upperA+\upperAlpha_0 =  \upperAlpha + \upperB + \upperJ\upperA \le (\upperB-1) + \upperB + \upperK\upperA = 2\upperB + \upperK \upperA-1,
	\end{equation*}
	and thus $\upperAlpha + \upperB + \upperJ\upperA \in \upperUU$, as desired.
	
	Second, we prove that $\upperWW_\upperN (\upperAlpha) \subset \upperUU$ if $0\le \upperN \le \upperNN$. Take an arbitrary element $\widehat\upperAlpha_\upperN + 2 \upperB + \upperJ\upperA \in \upperWW_\upperN (\upperAlpha)$. Recall that $\upperK + \upperTau_{\upperN} -\upperTau_{\upperN+1} \le \upperJ < \upperK$ and that $\upperTau_{\upperN+1} \le \upperK+1$. \Cref{t-for-main-sequence-lemma} implies that $\widehat\upperAlpha_\upperN + \upperB + \upperK' \upperA \in A$ for some $0 \le \upperK' < \upperTau_\upperN$. Note that $0\le (\upperK+1-\upperTau_{\upperN + 1})+(\upperTau_{\upperN}-\upperK'-1)= (\upperK + \upperTau_{\upperN} -\upperTau_{\upperN+1})-\upperK' \le \upperJ - \upperK' \le \upperK$. Hence, we have $(\widehat\upperAlpha_\upperN + 2 \upperB + \upperJ\upperA)-(\widehat\upperAlpha_\upperN + \upperB + \upperK' \upperA) = \upperB+(\upperJ - \upperK')\upperA \in \upperMM$, so $\widehat\upperAlpha_\upperN + 2 \upperB + \upperJ\upperA \notin \upperAA$. Besides, we have
	\begin{align*}
		&\widehat\upperAlpha_\upperN + 2 \upperB + \upperJ\upperA \le (\upperA-1) + 2 \upperB + (\upperK-1)\upperA = 2\upperB + \upperK \upperA-1, \\
		&\widehat\upperAlpha_\upperN + 2 \upperB + \upperJ\upperA \ge \widehat\upperAlpha_\upperN + 2 \upperB + (\upperK + \upperTau_{\upperN} -\upperTau_{\upperN+1})\upperA = \upperB+\upperK\upperA+(\upperB+(\upperTau_{\upperN}-\upperTau_{\upperN + 1})\upperA+\upperAlpha_{\upperN}),
	\end{align*}
	and \eqref{obs-N+1vsN} implies that the latter expression equals $\upperB + (\upperK+1) \upperA +\widehat\upperAlpha_{\upperN+1}$ if $\upperN<\upperNN$ and bounded from below by the same quantity if $\upperN=\upperNN$ due to a possible truncation of $\upperTau_{\upperNN + 1}$. In any case, we conclude that $\widehat\upperAlpha_\upperN + 2 \upperB + \upperJ\upperA \ge \upperB + (\upperK+1) \upperA$, and thus $\widehat\upperAlpha_\upperN + 2 \upperB + \upperJ\upperA \in \upperUU$, as desired.
\end{proof}

Next, let us study potential intersections between these sets.

\begin{Lemma}\label{alpha-sequences-are-disjoint}
	Given $\upperAlpha  = \upperEta_0 \upperA + \widehat \upperAlpha_0 \in \upperTT$, $\upperBeta  = \upperEta_0' \upperA + \widehat \upperBeta_0 \in \upperTT$ such that $0\le \upperEta_0, \upperEta_0' \le \upperK$, if $\widehat\upperAlpha_\upperN = \widehat\upperBeta_{\upperN'}$ for some $0 \le \upperN \le \upperNN$, $0 \le \upperN' \le \upperNN'$, then $\upperAlpha=\upperBeta$ and $\upperN = \upperN'$.
\end{Lemma}

\begin{proof}
	It is clear that if $\upperN \ge \upperN'$, then $\widehat\upperAlpha_{\upperN-\upperN'} = \widehat\upperBeta_{0}$, so we assume without loss of generality that $\upperN'=0$, i.e., that $\widehat \upperAlpha_\upperN = \widehat \upperBeta_0$. If $\upperN = 0$ and $\upperAlpha \neq \upperBeta$, we have  $|\upperAlpha - \upperBeta| = |\upperEta_0 - \upperEta_0'|\upperA \in \upperMM$ since $|\upperEta_0 - \upperEta_0'| \le \upperK$, a contradiction.
	
	If $\upperN > 0$, then $\widehat \upperAlpha_{\upperN-1} \not\in \upperII$ and \Cref{t-for-main-sequence-lemma} implies that $\widehat\upperAlpha_{\upperN - 1} + \upperB + \upperK' \upperA \in \upperAA$ for some $0 \le \upperK' < \upperTau_{\upperN-1}$. Note that $\upperBeta \in \upperTT$, and thus $\upperBeta < \upperB \le \widehat\upperAlpha_{\upperN - 1} + \upperB + \upperK' \upperA$. Moreover, observe that $ \upperTau_{\upperN} - \upperEta_0'-(\upperTau_{\upperN-1}-\upperK'-1) \le \upperK-0-0=\upperK$. Now it follows from \eqref{obs-N+1vsN} that
	$
		(\widehat\upperAlpha_{\upperN - 1} + \upperB + \upperK' \upperA) -\upperBeta = \upperB + \upperK' \upperA - \upperEta_0' \upperA + (\widehat\upperAlpha_{\upperN - 1} - \widehat\upperAlpha_{\upperN}) = (\upperTau_{\upperN} - \upperEta_0'-\upperTau_{\upperN-1}+\upperK'+1)\upperA \in \upperMM,
	$ which is a contradiction again.
\end{proof}

\begin{Lemma}\label{injective-correspondence}
	For all $\upperAlpha, \upperBeta \in \upperTT$, the following statements are valid:
	\begin{itemize}
		\item $\upperVV(\upperAlpha) \cap \upperVV(\upperBeta) = \varnothing$ if $\upperAlpha \neq \upperBeta$.
		\item $\upperWW_\upperN(\upperAlpha) \cap \upperVV(\upperBeta) = \varnothing$ for all $0 \le \upperN \le \upperNN$.
		\item $\upperWW_\upperN (\upperAlpha) \cap \upperWW_{\upperN'} (\upperBeta)=\varnothing$ for all $0 \le \upperN \le \upperNN$, $0 \le \upperN' \le \upperNN'$ unless $\upperAlpha = \upperBeta$ and $\upperN = \upperN'$.
	\end{itemize}
\end{Lemma}
\begin{proof}
    Let us prove these statements one by one: 
    \begin{enumerate}
        \item If $\upperAlpha \neq \upperBeta$ and $\upperVV(\upperAlpha) \cap \upperVV(\upperBeta) \neq \varnothing$, then $\upperAlpha + \upperB + \upperJ\upperA = \upperBeta + \upperB + \upperJ'\upperA$ for some $0 \le \upperJ, \upperJ' \le \upperK$, and thus $|\upperAlpha - \upperBeta| = |\upperJ - \upperJ'|\upperA \in \upperMM$, a contradiction.
        \item Suppose that $\upperWW_\upperN (\upperAlpha) \cap \upperVV (\upperBeta) \neq \varnothing$, i.e., that $\widehat\upperAlpha_\upperN + 2\upperB + \upperJ\upperA = \upperBeta + \upperB + \upperJ' \upperA$ for some $0 \le \upperJ' \le \upperK$ and $\upperK+\upperTau_{\upperN}-\upperTau_{\upperN + 1} \le \upperJ \le \upperK$. \Cref{t-for-main-sequence-lemma} yields that $\widehat\upperAlpha_\upperN + \upperB + \upperK' \upperA \in \upperAA$ for some $0 \le \upperK' < \upperTau_\upperN$. Note that $\upperBeta \in \upperTT$, and thus $\upperBeta < \upperB \le \widehat \upperAlpha_\upperN + \upperB + \upperK'\upperA$. Moreover, we have $\upperK'+\upperJ' - \upperJ \le (\upperTau_{\upperN}-1)+\upperK - (\upperK+\upperTau_{\upperN} -\upperTau_{\upperN + 1}) = \upperTau_{\upperN + 1}-1 \le \upperK$ and thus $(\widehat\upperAlpha_\upperN + \upperB + \upperK' \upperA) - \upperBeta = (\upperK'+\upperJ' - \upperJ)\upperA \in \upperMM$, a contradiction.
        \item Suppose that $\upperWW_\upperN (\upperAlpha) \cap \upperWW_{\upperN'} (\upperBeta) \neq \varnothing$, i.e., that $\widehat\upperAlpha_\upperN + 2\upperB + \upperJ\upperA = \widehat\upperBeta_{\upperN'} + 2\upperB + \upperJ'\upperA$ for some $\upperJ, \upperJ'$. Since $|\widehat\upperAlpha_\upperN - \widehat\upperBeta_{\upperN'}|< \upperA$, we conclude that $\upperJ =  \upperJ'$ and thus $\widehat\upperAlpha_\upperN = \widehat\upperBeta_{\upperN'}$. Now \Cref{alpha-sequences-are-disjoint} completes the proof. \qedhere
    \end{enumerate}
\end{proof}

It easily follows from the last two lemmas that the images of distinct $\upperAlpha, \upperBeta \in \upperTT$ are disjoint, whether they are elements of $\upperII$ or subsets of $\upperUU$, respectively. Hence, it remains only to check that if $\upperAlpha = \upperEta_0 \upperA + \widehat \upperAlpha_0 \in \upperTT$ is mapped into $\upperUU (\upperAlpha) \subset \upperUU$, then $|\upperUU (\upperAlpha)| = \upperK+1$. If $\upperEta_0 > \upperK$, then there is nothing to do. Otherwise, we have $\upperUU(\upperAlpha) = \upperVV (\upperAlpha) \cup \bigcup_{\upperN = 0}^{\upperNN} \upperWW_\upperN(\upperAlpha)$, and this union is disjoint by \Cref{injective-correspondence}. Thus 
\begin{equation*}
    \left| \upperUU (\upperAlpha) \right| = \left| \upperVV (\upperAlpha) \right| + \sum_{\upperN = 0}^{\upperNN} \left| \upperWW_\upperN(\upperAlpha) \right| = \upperEta_0 + \sum_{\upperN = 0}^{\upperNN} (\upperTau_{\upperN + 1}-\upperTau_{\upperN}) = \upperEta_{\upperNN+1}= \upperK + 1. 
\end{equation*}

\section{Concluding remarks} \label{SecConc}

Though the tightness of \Cref{main-conjecture} is the main open problem of our paper, we briefly discuss some more.

\vspace{1mm}
\noindent \textbf{Coverings.}
For a finite $S \subset \Z$, its \textit{covering density} $d_c(S)$ is defined as the minimum lower density of $A\subset \Z$ such that $\bigcup_{a \in A} (a+S)=\Z$. Since Newman~\cite{newman1967complements} introduced this notion in 1967, it has been extensively studied, often under various different terms, see \cite{axe2019polychromatic,bollobas2011covering,hainzl2022finding,huang2019domination,schmidt2008covering,schmidt2010covering}. It would be interesting to determine the covering density of $S = \{0, \theora, \dots, \theork\theora, \theork\theora + \theorb, \dots, \theork\theora + \theorm\theorb\}$ as well. In case $k=m=1$, this problem goes back to Schmidt and Tuller~\cite{schmidt2008covering} and was recently solved in~\cite{frankl2023solution}, while the case $k=b=1$ was handled in~\cite{huang2022domination}. 

\vspace{1mm}
\noindent \textbf{Symmetric $4$-element sets.}
In 2004, Liu and Zhu~\cite[Conjecture~4.1]{liuzhu2004} stated a conjecture on the packing density of $S=\{0,a,b,a+b\}$, which remains open. Here we present its covering counterpart, which is also connected to a certain problem in Ramsey theory, see~\cite[Section~2.1]{frankl2024max} and~\cite{frankl2023solution}.
\begin{Conjecture}
	Let $a<b$ be coprime positive integers. Then for $S=\{0,a,b,a+b\}$, we have
	\begin{equation*}
		d_p(S) =
		\begin{dcases}
			\frac{1}{4}, & b-a \mbox{ is odd}, \\
			\frac{\lfloor(ab+b-a)/4\rfloor}{ab+b-a}, & b-a \mbox{ is even};
		\end{dcases}
		\hspace{7mm} d_c(S) =
		\begin{dcases}
			\frac{1}{4}, & b-a \mbox{ is odd}, \\
			\frac{\lceil(ab+a+b)/4\rceil}{ab+a+b}, & b-a \mbox{ is even}.
		\end{dcases}
	\end{equation*}
\end{Conjecture}
\noindent
In fact, if $b-a$ is odd, then it is easy to see that $S=\{0,a,b,a+b\}$ tiles $\Z$ and thus $d_p(S)=d_s(S)=1/4$. We state these equalities as parts of the conjecture only for completeness. Moreover, if $b-a$ is even, then for both packing and covering parts of this conjecture, it is not very difficult to find an explicit construction of the desired density. Proving its optimality seems much more challenging.

\vspace{1mm}
\noindent \textbf{Kappa values.} For every finite $M \subset \N$, its $\mu(M)$ is bounded from below by a certain parameter $\kappa(M)$, related to the `lonely runner conjecture', and the equality $\mu(M)=\kappa(M)$ often holds, see~\cite{liu2008rainbow} and the references therein. In particular, it holds for $M=\{a,b,a+b\}$, see \cite[Theorem~5.1]{liuzhu2004}. We wonder if this is the case in a more general setting where $M=\{\loweri\lowera+\lowerj\lowerb: 0 \le \loweri \le \lowerk,\, 0 \le \lowerj \le \lowerm,\, \loweri+\lowerj>0\}$ as well. We only note that if $a-b$ and $k+m+1$ are coprime, then our proof of \Cref{main-conjecture} also provides the same lower bound on $\kappa(M)$ as it does on $\mu(M)$.

%\subsection*{Acknowledgments}

\vspace{3mm}
\noindent
{\large \bf Acknowledgments.} Arsenii Sagdeev was supported by ERC Advanced Grant `GeoScape' No. 882971.

\printbibliography
\end{document}